\documentclass[11pt]{article}
\usepackage{amsthm,amsfonts,amsmath,amssymb}
\usepackage{mdwlist}
\usepackage{array,color}
\newtheorem*{thm A}{Theorem A}
\newtheorem*{thm B}{Theorem B}
\newtheorem*{thm C}{Theorem C}
\newtheorem*{prob}{Open Problems}

\newtheorem{thm}{Theorem}[section]

\newtheorem{lem}[thm]{Lemma}
\newtheorem{cor}[thm]{Corollary}

\newtheorem{exam}[thm]{Example}

\def\demo{\noindent{\bf Proof}\hskip10pt}

\def\qed{\hfill $\Box$}
\def\lg{\langle}
\def\rg{\rangle}

\begin{document}
\title{Finite groups with the most Chermak-Delgado measures of subgroups\thanks{This work was supported by NSFC (No. 12371022 \&12271318)}}
\author{Guojie Liu\ \ Haipeng Qu\ \ Lijian An\thanks{Corresponding author. e-mail: anlj@sxnu.edu.cn}\\
Shanxi Key Laboratory of Cryptography and Data Security, Shanxi Normal University\\
Taiyuan, Shanxi 030031, P. R. China\\
 }

\maketitle

\begin{abstract}

Let $G$ be a finite group and $H\leq G$. The Chermak-Delgado measure of $H$ is defined as the number $|H|\cdot|C_{G}(H)|$.
In this paper, we identify finite groups
that exhibit the  maximum number of Chermak-Delgado measures under some specific conditions.

\medskip


\end{abstract}

\baselineskip=16pt

\section{Introduction}
Throughout this paper, let $G$ be a finite group and $\mathcal{L}(G)$ be the subgroup lattice
of $G$.  $m_G$ is a map from $\mathcal{L}(G)$ to $\mathbb{N}^{*}$, given by $$m_G(H)=|H|\cdot |C_G(H)|,$$ where $H\in \mathcal{L}(G)$. Following Isaacs \cite{I}, $m_G(H)$ is called the Chermak-Delgado measure of $H$ (in $G$), and the subset of $\mathcal{L}(G)$ composed of subgroups with the maximum Chermak-Delgado measure is called the Chermak-Delgado lattice of $G$. In recent years, there has been a growing interest in understanding this lattice (see e.g. [1-5,~7-13,~15,~18-24]).

Following T${\rm \breve{a}}$rn${\rm \breve{a}}$uceanu, we use $|\mathrm{Im}(m_{G})|$ to denote the number of distinct Chermak-Delgado measures of subgroups in $G$.
T${\rm \breve{a}}$rn${\rm \breve{a}}$uceanu \cite{T2} points out that $|\mathrm{Im}(m_{G})|\geqslant 2$
for every nontrivial finite group $G$, and provides some properties of $G$ with $|\mathrm{Im}(m_{G})|=2$.

It is natural to consider the upper bound of $|\mathrm{Im}(m_{G})|$. Notice that if $G$ is an abelian group of order $p^{n-1}$, then $|\mathrm{Im}(m_{G})|=n$. Therefore, for finite groups, $|\mathrm{Im}(m_{G})|$ has no upper bound.
In this paper, we study the upper bound of $|\mathrm{Im}(m_{G})|$ under some specific conditions. For convenience, we give the following notations:
$$\mathrm{Im}_{\rm max}(n)={\rm max}\{|\mathrm{Im}(m_{G})|\mid G ~{\rm is~a~group~of~order} ~n\};$$
$$\mathrm{Im}_{\rm max}(n,\mathcal{N})={\rm max}\{|\mathrm{Im}(m_{G})|\mid  G ~{\rm is~a~nilpotent ~group~of~order} ~n \}.$$

First, we consider the case $n=p^k$, where $p$ is a prime and $k$ is a positive integer. Following result is obtained:

\begin{thm A}
Let $n=p^k$, where $p$ is a prime and $k$ is a positive integer. Suppose that $G$ is a group of order $n$. Then

$(1)$ If $k<5$, then $\mathrm{Im}_{\rm max}(n)=k+1$. Moreover, $|\mathrm{Im}(m_{G})|=\mathrm{Im}_{\rm max}(n)$ if and only if $G$ is abelian.

$(2)$ If $k>5$, then $\mathrm{Im}_{\rm max}(n)=2k-4$. Moreover, $|\mathrm{Im}(m_{G})|=\mathrm{Im}_{\rm max}(n)$ if and only if $G$ is of maximal class with an abelian subgroup of index $p$ and a uniform element of order $p$.

$(3)$ If $k=5$, then $\mathrm{Im}_{\rm max}(n)=6$. Moreover, $|\mathrm{Im}(m_{G})|=\mathrm{Im}_{\rm max}(n)$ if and only if $G$ is abelian or $G$ is of maximal class with an abelian subgroup of index $p$ and a uniform element of order $p$.

\end{thm A}

Using above result, following result on $\mathrm{Im}_{\rm max}(n,\mathcal{N})$ is also obtained:

\begin{thm B}
Let $n=p_{1}^{k_{1}}p_{2}^{k_{2}}\cdots p_{s}^{k_{s}}$, where $p_1,\ldots,p_s$ are distinct prime factors of $n$, $k_{i}\in \mathbb{N}^{*},~i=1,2,\ldots,s$. Then

$(1)$ $\mathrm{Im}_{\rm max}(n,\mathcal{N})=\prod\limits_{k_{i}\leqslant 4}(k_{i}+1)\cdot\prod\limits_{k_{i}\geqslant 5}(2k_{i}-4),i=1,2,\ldots,s;$

$(2)$ Let $G$ be~a~nilpotent group~of~order $n$. Then $|\mathrm{Im}(m_{G})|=\mathrm{Im}_{\rm max}(n,\mathcal{N})$ if and only if $G=P_{1}\times P_{2}\times\cdots\times P_{s}$, where $P_{i}\in\mathrm{Syl}_{p_{i}}(G)$, $P_{i}$ is abelian when $k_{i}<5$, $P_{i}$ is of maximal class with an abelian subgroup of index $p_{i}$ and a uniform element of order $p_{i}$ when $k_{i}>5$, and $P_{i}$ is abelian or $P_{i}$ is of maximal class with an abelian subgroup of index $p_{i}$ and a uniform element of order $p_{i}$ when $k_{i}=5$, $i=1,\ldots,s$.

\end{thm B}

We also consider the case where $n$ is square-free. Following result is obtained:

\begin{thm C}
Let $n=p_{1}p_{2}\cdots p_{s}$, where $p_1,\ldots,p_s$ are distinct prime factors of $n$. Suppose that $G$ is a group of order $n$. Then $\mathrm{Im}_{\rm max}(n)=2^s.$ Moreover, $|\mathrm{Im}(m_{G})|=\mathrm{Im}_{\rm max}(n)$ if and only if $G$ is abelian.

\end{thm C}

\section{Preliminaries}
Let $G$ be a finite group. We use $c(G)$ to denote the nilpotency class of $G$. For a nilpotent group $G$, let
$$G=K_1(G)> G'=K_2(G)> K_3(G)> \cdots> K_{c+1}(G)=1$$
be the lower central series of $G$, where $c=c(G)$ and $K_{i+1}(G)=[K_i(G),G]$ for $1\leqslant i\leqslant c$.

A non-abelian group $G$ of order $p^n$ is said to be maximal class if $c(G)=n-1$. We say
that $s\in G$ is a uniform element if $s\notin \bigcup\limits_{i=2}^{n-2} C_{G}(K_{i}(G)/K_{i+2}(G))$.

\begin{lem}\label{mc1}{\rm \cite[Proposition 1.8]{YB}}
Let $G$ be a non-abelian $p$-group. If $A<G$ of order $p^2$ is such that $C_{G}(A)=A$, then $G$ is of maximal class.
\end{lem}

\begin{lem}\label{mc2}{\rm \cite[Theorem 3.15(ii)]{GA}}
Let $G$ be a $p$-group of maximal class and suppose that $G$ has a uniform element $s$. Then $|C_{G}(s)|=p^2$.
\end{lem}

\begin{lem}{\rm \cite[VI, 1.8 Hauptsatz]{HB}(P. Hall)}\label{s1}
Suppose that $G$ is a finite solvable group, and let $\pi$ be an arbitrary set of primes. Then $G$ has a Hall $\pi$-subgroup, any two Hall $\pi$-subgroups of $G$ are conjugate in $G$, and any $\pi$-subgroup of $G$ is contained in a Hall $\pi$-subgroup of $G$.
\end{lem}

\begin{lem}{\rm \cite[Corollary 5.15]{I}}\label{s3}
Let $G$ be a finite group, and suppose that all Sylow subgroups of $G$ are cyclic. Then $G$ is solvable.
\end{lem}

\section{Proof of Main Theorem}
Before we prove the main theorem, we will first provide some conclusions that will be employed in the following proofs.

\begin{thm}\label{thm1}
Let $G=H\times K$, $A\leq H$ and $B\leq K$. Then $m_{G}(A\times B)=m_{H}(A)\cdot m_{K}(B)$.
\end{thm}
\demo
We have $$C_{G}(A\times B)=C_{G}(A)\cap C_{G}(B)=(C_{H}(A)\times K)\cap (H\times C_{K}(B))=C_{H}(A)\times C_{K}(B).$$
 It follows that
\begin{equation*}
 \begin{aligned}
 m_{G}(A\times B)=&|A\times B|\cdot|C_{G}(A\times B)|=|A|\cdot|B|\cdot|C_{H}(A)\times C_{K}(B)|\\
                 =&|A|\cdot|C_{H}(A)|\cdot|B|\cdot|C_{K}(B)|=m_{H}(A)\cdot m_{K}(B).
 \end{aligned}
\end{equation*}\qed

\begin{cor}\label{cor2}
Let $G=H\times K$. Then for all $a\in \mathrm{Im}(m_{H})$ and $b\in \mathrm{Im}(m_{K})$, we have $ab\in \mathrm{Im}(m_{G}).$
\end{cor}

\demo
For all $a\in \mathrm{Im}(m_{H})$ and $b\in \mathrm{Im}(m_{K})$, there exist $A\leq H$ and $B\leq K$ such that $m_{H}(A)=a$ and $m_{K}(B)=b$. By Theorem \ref{thm1}, we have $ab=m_{G}(A\times B)\in \mathrm{Im}(m_{G}).$\qed

\begin{lem}\label{dr}
Let $G=H\times K$. If $(|H|,|K|)=1$, then for any $M\leq G$, there exist $A\leq H$ and $B\leq K$ such that $M=A\times B$.
\end{lem}
\demo
We assert that $M=(H\cap M)\times (K\cap M)$. Obviously, we have $M\supseteq(H\cap M)\times (K\cap M)$. Therefore, we only need to prove $M\subseteq(H\cap M)\times (K\cap M)$. For any $g\in M$, there exist $h_1\in H$ and $k_1\in K$ such that $g=h_1k_1$. Let $o(h_{1})=l,~o(k_{1})=t$. Since $(|H|,|K|)=1$, $(l,t)=1$. Then there exist $u,v\in \mathbb{Z}$ such that $g=g^{ul+vt}={h_1}^{vt}{k_1}^{ul}\in \langle h_{1}^{t}, k_{1}^{l}\rangle=\langle h_{1}\rangle\times\langle k_{1}\rangle=\langle g^t\rangle\times\langle g^l\rangle\subseteq(H\cap M)\times (K\cap M)$. It follows that $M=(H\cap M)\times (K\cap M)$. That is $M=A\times B$, where $A=H\cap M$ and $B=K\cap M$.\qed

\begin{lem}\label{lem1}
Let $G=H\times K$. Then

$(1)$ $|\mathrm{Im}(m_{G})|\geqslant|\mathrm{Im}(m_{H})|+|\mathrm{Im}(m_{K})|-1;$

$(2)$ If $(|H|,|K|)=1$, then $|\mathrm{Im}(m_{G})|=|\mathrm{Im}(m_{H})|\cdot|\mathrm{Im}(m_{K})|$.

\end{lem}

\demo
Let $$\mathrm{Im}(m_{H})=\{a_1,a_2,\cdots,a_s\}~\mbox{and}~\mathrm{Im}(m_{K})=\{b_1,b_2,\cdots,b_t\},$$ where $a_1<a_2<\cdots<a_s~\mbox{and}~b_1<b_2<\cdots<b_t$. By Corollary \ref{cor2}, $$a_{i}b_j\in \mathrm{Im}(m_{G}),~\mbox{for}~ i\in \{1,2,\ldots,s\}~\mbox{and}~j\in \{1,2,\ldots,t\}.$$

$(1)$  Since $a_1b_1<a_1b_2<\cdots<a_1b_t<a_2b_t<\cdots<a_sb_t,$ $$\{a_1b_1,a_1b_2,\ldots,a_1b_t,a_2b_t,\ldots,a_sb_t\}\subseteq \mathrm{Im}(m_{G}).$$ Thus $$|\mathrm{Im}(m_{G})|\geqslant|\mathrm{Im}(m_{H})|+|\mathrm{Im}(m_{K})|-1.$$

$(2)$  Let $M\leq G$. By Lemma \ref{dr}, there exist $A\leq H$ and $B\leq K$ such that $M=A\times B$. Therefore, $m_{G}(M)=m_{H}(A)\cdot m_{K}(B)$ by Theorem \ref{thm1}. It follows that $m_{G}(M)=a_ib_j$ for some $i$ and $j$. On the other hand, since $(|H|,|K|)=1$, $(a_{i},b_{j})=1$ where $1\leqslant i\leqslant s$ and $1\leqslant j\leqslant t$. Hence $a_{i}b_{j}\neq a_{i'}b_{j'}$ for $i\neq i'$ or $j\neq j'$, where $1\leqslant i,i'\leqslant s$ and $1\leqslant j,j'\leqslant t$.
 Thus $\mathrm{Im}(m_{G})=\{a_ib_j\mid 1\leqslant i\leqslant s, 1\leqslant j\leqslant t\}$ and $|\mathrm{Im}(m_{G})|=st=|\mathrm{Im}(m_{H})|\cdot|\mathrm{Im}(m_{K})|$.

\qed

\begin{lem}\label{ge}
Let $G$ be a finite group.  Suppose that $H\leq G$. Then $m_{G}(H)=m_{G}(H^{g})$ for all $g\in G$.
\end{lem}

\demo
Since $C_{G}(H^{g})=C_{G^{g}}(H^{g})=C_{G}(H)^{g}$, $|C_{G}(H^{g})|=|C_{G}(H)^{g}|=|C_{G}(H)|$.
Hence $m_{G}(H^{g})=|H^{g}|\cdot|C_{G}(H^{g})|=|H|\cdot|C_{G}(H)|=m_{G}(H)$.\qed

\begin{cor}\label{ge1}
Let $G$ be a finite non-abelian group. Then $|\mathrm{Im}(m_{G})|\leqslant c-1$, where $c$ denotes the total number of conjugacy classes of subgroups in $G$.
\end{cor}
\demo
Since $G$ is non-abelian, $G$ and $Z(G)$ are not conjugate.
Since $m_{G}(G)=m_{G}(Z(G))$, by Lemma \ref{ge}, $|\mathrm{Im}(m_{G})|\leqslant c-1$.

\begin{lem}\label{a1}
Let $G$ be an abelian group of order $p^{k}$. Then $|\mathrm{Im}(m_{G})|=k+1$.
\end{lem}
\demo
Let $H\leq G$ such that $|H|=p^s$, where $0\leqslant s\leqslant k$. Since $G$ is abelian, $C_{G}(H)=G$. Hence $$m_{G}(H)=|H|\cdot|C_{G}(H)|=|H|\cdot|G|=p^{k+s}.$$
It follows that $|\mathrm{Im}(m_{G})|=k+1$. \qed

\begin{lem}\label{na1}
Let $G$ be a non-abelian group of order $p^{3}$. Then $|\mathrm{Im}(m_{G})|=2$.
\end{lem}
\demo
Since $G$ is a non-abelian group of order $p^{3}$, $|Z(G)|=p$. Hence $m_{G}(G)=m_{G}(Z(G))=m_{G}(H)=p^4$, where $H\leq G$ and $|H|=p^2$. For any $K\leq G$ such that $|K|=p$ and $K\neq Z(G)$, we have $|C_{G}(K)|=p^2$. Thus $m_{G}(K)=m_{G}(1)=p^3$. Hence $|\mathrm{Im}(m_{G})|=2$. \qed

\begin{lem}\label{na2}
Let $G$ be a non-abelian group of order $p^{k}(k>3)$. Then

 $(1)$ Let $H\leq G$. Then $m_{G}(H)\geqslant p^{3}$, where ``=" holds if and only if $|H|=|Z(G)|=p$ and $C_G(H)=HZ(G)$.

 $(2)$ Let $H\leq G$. Then $m_{G}(H)\leqslant p^{2k-2}$, where ``=" holds if and only if $G$ has an abelian maximal subgroup.

 $(3)$ $|\mathrm{Im}(m_{G})|=2k-4$ if and only if $G$ is of maximal class with an abelian maximal subgroup and a uniform element $x$ of order $p$.
\end{lem}
\demo
$(1)$ If $H\leq Z(G)$, then $C_G(H)=G$ and $$m_{G}(H)=|H|\cdot|C_{G}(H)|=|H|\cdot|G|\geqslant p^k>p^3.$$
If $H\nleq Z(G)$ and $H$ is non-abelian, then $|H|\geqslant p^3$ and $$m_{G}(H)= |H|\cdot|C_{G}(H)|\geqslant|H|\cdot|Z(G)|>p^3.$$ If $H\nleq Z(G)$ and $H$ is abelian, then $|H|\geqslant p$, $C_G(H)\ge HZ(G)>H$ and $$m_{G}(H)= |H|\cdot|C_{G}(H)|\geqslant|H|\cdot|HZ(G)|\geqslant p^3.$$ In conclusion, $m_{G}(H)\geqslant p^{3}$, where ``=" holds if and only if $|H|=|Z(G)|=p$ and $C_G(H)=HZ(G)$.

$(2)$ First, consider the case where $G$ has an abelian maximal subgroup $A$. Then $m_{G}(A)=|A|\cdot|C_{G}(A)|=|A|\cdot|A|=p^{2k-2}.$ In the following, we prove that $m_G(H)\le p^{2k-2}$ for $H\le G$.
 Since $G$ is non-abelian, $|Z(G)|\leqslant p^{k-2}$. Hence $$m_{G}(G)=|G|\cdot|Z(G)|\leqslant p^k\cdot p^{k-2}=p^{2k-2}.$$ If $|H|=p^{k-1}$, then $$m_{G}(H)=|H|\cdot|C_{G}(H)|\leqslant p^{k-1}\cdot p^{k-1}=p^{2k-2}.$$  If $|H|\leqslant p^{k-2}$, then $$m_{G}(H)=|H|\cdot|C_{G}(H)|\leqslant p^{k-2}\cdot p^{k}=p^{2k-2}.$$
Next, consider the case where $G$ has no abelian maximal subgroup. In this case, $|Z(G)|\leqslant p^{k-3}$. Hence $$m_{G}(G)=|G|\cdot|Z(G)|\leqslant p^k\cdot p^{k-3}=p^{2k-3}.$$
If $H\le Z(G)$, then $|H|\leqslant p^{k-3}$ and $$m_{G}(H)=|H|\cdot|C_{G}(H)|=|H|\cdot|G|\leqslant p^{k-3}\cdot p^{k}=p^{2k-3}.$$
If $H\nleq Z(G)$ and $H$ is a proper subgroup of $G$, then $|H|\leqslant p^{k-1}$, $|C_{G}(H)|\leqslant p^{k-1}$. Particularly, if $|H|=p^{k-1}$, then $|C_{G}(H)|<p^{k-1}$ since $G$ has no abelian maximal subgroup. Hence we have $$m_{G}(H)=|H|\cdot|C_{G}(H)|< p^{k-1}\cdot p^{k-1}=p^{2k-2}.$$
In conclusion, $m_{G}(H)\leqslant p^{2k-2}$, where ``=" holds if and only if $G$ has an abelian maximal subgroup.

$(3)$ ($\Rightarrow$) By $(1)$ and $(2)$, we have  $$\mathrm{Im}(m_{G})\subseteq\{p^3,p^4,\ldots,p^{2k-2}\}.$$ Since $|\mathrm{Im}(m_{G})|=2k-4$, $\mathrm{Im}(m_{G})=\{p^3,p^4,\ldots,p^{2k-2}\}.$ By (2), there exists $G_{1}\leq G$ such that $|G_{1}|=p^{k-1}$ and $C_{G}(G_{1})=G_{1}$. By (1), there exists $H\le G$ such that $|H|=|Z(G)|=p$ and $C_G(H)=HZ(G)$. Let $C_G(H)=K$. Then $|K|=p^2$ and $C_G(K)=C_G(H)=K$. By Lemma \ref{mc1}, $G$ is of maximal class. Hence $C_{G}(K_{i}(G)/K_{i+2}(G))<G$ for $i=2,3\dots,k-2$.
Since $K_{i}(G)< G_{1}$ and $G_1$ is abelian, $G_{1}\leq C_{G}(K_{i}(G)/K_{i+2}(G))$. It follows that $G_{1}=C_{G}(K_{i}(G)/K_{i+2}(G))$, and hence $$G_{1}=\bigcup\limits_{i=2}^{k-2} C_{G}(K_{i}(G)/K_{i+2}(G)).$$  Let $H=\lg x\rg$. Since $C_G(H)=HZ(G)$, $x\not\in G_1$. That is, $x$ is a uniform element.

($\Leftarrow$) By Lemma \ref{mc2}, $|C_{G}(x)|=p^2$. Hence $m_{G}(\langle x\rangle)=|\langle x\rangle|\cdot|C_{G}(x)|=p^3$. Let $H=\langle x\rangle K_{i}(G)$, where $i=2,3,\ldots,k-2$. Then $|H|=p^{k-i+1}$ and $$C_{G}(H)=C_{G}(\langle x\rangle K_{i}(G))=C_{G}(x)\cap C_{G}(K_{i}(G))=\langle x\rangle Z(G)\cap C_{G}(K_{i}(G))=Z(G).$$
Hence $$m_{G}(H)=|H|\cdot|C_{G}(H)|=|H|\cdot|Z(G)|=p^{k-i+1}\cdot p=p^{k-i+2},~2\leqslant i\leqslant k-2.$$ Let $G_1$ be the abelian maximal subgroup of $G$ and $T$ be a subgroup of $G_1$ of order $p^t~,2\leqslant t\leqslant k-1$. Then $C_{G}(T)=G_{1}$. If not, then $C_{G}(T)=G$. It follows that $T\leq Z(G)$, this is a contradiction. Hence $$m_{G}(T)=|T|\cdot|C_{G}(T)|=p^{t+k-1},~2\leqslant t\leqslant k-1.$$ Then $\mathrm{Im}(m_{G})=\{p^3,p^4,\ldots,p^{2k-2}\}$. Hence $|\mathrm{Im}(m_{G})|=2k-4$. \qed

\bigskip
\noindent {\bf Proof of Theorem A.} By Lemma \ref{a1}, Lemma \ref{na1} and Lemma \ref{na2}, if $k<5$, then $\mathrm{Im}_{\rm max}(n)=\max\{k+1, 2k-4\}=k+1$ and $|\mathrm{Im}(m_{G})|=k+1$ if and only if $G$ is abelian;
if $k>5$, then $\mathrm{Im}_{\rm max}(n)=\max\{k+1, 2k-4\}=2k-4$ and $|\mathrm{Im}(m_{G})|=2k-4$ if and only if $G$ is of maximal class with an abelian maximal subgroup and a uniform element of order $p$;
if $k=5$, then $\mathrm{Im}_{\rm max}(n)=\max\{k+1, 2k-4\}=6$ and $|\mathrm{Im}(m_{G})|=6$ if and only if $G$ is abelian or $G$ is of maximal class with an abelian maximal subgroup and a uniform element of order $p$.\qed

\bigskip

\noindent {\bf Proof of Theorem B.} For any nilpotent group $G$ of order $n$, we have $G=P_{1}\times P_{2}\times\cdots\times P_{s}$, where $P_{i}\in\mathrm{Syl}_{p_{i}}(G),~i=1,2,\ldots,s$. Since $(|P_{m}|,|P_{n}|)=1$, where $m,n\in \{1,2,\ldots,s\}$ and $m\neq n$,  $|\mathrm{Im}(m_{G})|=\prod\limits_{i=1}^{s}|\mathrm{Im}(m_{P_{i}})|$ by Lemma \ref{lem1}(2). Thus, we have $$|\mathrm{Im}(m_{G})|=\mathrm{Im}_{\rm max}(n,\mathcal{N})\Longleftrightarrow |\mathrm{Im}(m_{P_{i}})|=\mathrm{Im}_{\rm max}(p_{i}^{k_{i}},\mathcal{N}),~i=1,2,\ldots,s.$$

 $(1)$  Hence $$\mathrm{Im}_{\rm max}(n,\mathcal{N})=\prod\limits_{i=1}^{s}\mathrm{Im}_{\rm max}(p_{i}^{k_{i}},\mathcal{N}).$$ By Theorem A$(1)$, $$\mathrm{Im}_{\rm max}(n,\mathcal{N})=\prod\limits_{k_{i}\leqslant 4}(k_{i}+1)\cdot\prod\limits_{k_{i}\geqslant 5}(2k_{i}-4),i=1,2,\ldots,s.$$

$(2)$ By Theorem A, we have $|\mathrm{Im}(m_{G})|=\mathrm{Im}_{\rm max}(n,\mathcal{N})$ if and only if $$G=P_{1}\times P_{2}\times\cdots\times P_{s},$$ where $P_{i}\in\mathrm{Syl}_{p_{i}}(G)$, $P_{i}$ is abelian when $k_{i}<5$, $P_{i}$ is of maximal class with an abelian subgroup of index $p_{i}$ and a uniform element of order $p_{i}$ when $k_{i}>5$, and $P_{i}$ is abelian or $P_{i}$ is of maximal class with an abelian subgroup of index $p_{i}$ and a uniform element of order $p_{i}$ when $k_{i}=5$, $i=1,\ldots,s$. \qed

\bigskip
\noindent {\bf Proof of Theorem C.} Since all Sylow subgroups of $G$ are cyclic, $G$ is solvable by Lemma \ref{s3}. It follows from Lemma \ref{s1} that all Hall $\pi$-subgroups of $G$ are conjugate, where $\pi$ is an arbitrary set of prime factors of $n$. Since any subgroup of $G$ is Hall subgroup of $G$, all subgroups of the same order of $G$ are conjugate. Therefore the Chermak-Delgado measures of the same order subgroups of $G$ are equal by Lemma \ref{ge}. Since the number of subgroup orders of $G$ is $2^s$ by Lemma \ref{s1}, $\mathrm{Im}_{\rm max}(n)\leqslant 2^s$. If $G$ is abelian, then for any $H\leq G$ we have $C_{G}(H)=G$. Hence the Chermak-Delgado measure of each order subgroup of $G$ is distinct. Thus $|\mathrm{Im}(m_{G})|$ is equal to the number of subgroup orders of $G$. That is $|\mathrm{Im}(m_{G})|=2^s$. If $G$ is non-abelian, then $|Z(G)|<|G|$. Moreover, $|\mathrm{Im}(m_{G})|<2^s$ since $m_{G}(G)=m_{G}(Z(G))$. In conclusion, $\mathrm{Im}_{\rm max}(n)=2^s$ and $|\mathrm{Im}(m_{G})|=2^s$ if and only if $G$ is abelian.\qed

\section{The Upper Bound of $|\mathrm{Im}(m_{G})|$}
After establishing the exact upper bound for $|\mathrm{Im}(m_{G})|$ of group $G$ under certain conditions, we now further investigate the upper bound of $|\mathrm{Im}(m_{G})|$ in order to later determine the exact upper bound for more general groups $G$. We use $\tau(n)$ to denote the number of factors of positive integer $n$.

\begin{thm}
Let $G$ be a group of order $n$. Suppose that $|Z(G)|=m$. Then $|\mathrm{Im}(m_{G})|\leqslant \frac{(\tau(n)-1)(\tau(n)-2)}{2}+\tau(m)$.
\end{thm}
\demo
Let $H\leq G$. Since $|H|\mid n$ and $|C_{G}(H)|\mid n$, $m_{G}(H)=|H|\cdot|C_{G}(H)|$ is the product of two factors of $n$. In particular, for any $H\leq Z(G)$, we have $m_{G}(H)=|H|\cdot|G|=n|H|$. Thus $|\{m_{G}(H)\mid H\le Z(G)\}|=\tau(m)$.
Excluding the cases of $1,n$, there are $\tau(n)-2$ non-trivial factors of $n$. The number of ways to choose two distinct non-trivial factors to form a product is$\binom{\tau(n)-2}{2}$, and the number of ways to choose a non-trivial factor to form a square product is \(\tau(n)-2\). Adding the values $m_G(H)(H\le Z(G))$ to the count, we have $$|\mathrm{Im}(m_G)|\leqslant\binom{\tau(n)-2}{2}+(\tau(n)-2)+\tau(m)
=\frac{(\tau(n)-1)(\tau(n)-2)}{2}+\tau(m).$$ \qed

\begin{lem}\label{t2}
Let $G=H\times K$, where $K$ is abelian. Then $|\mathrm{Im}(m_{G})|\leqslant|\mathrm{Im}(m_{H})|\cdot|\mathrm{Im}(m_{K})|$.

\end{lem}
\demo Let $A\leq G$. By modular law, $AK=(AK\cap H)\times K$. Hence $$\frac{|A||K|}{|A\cap K|}=|AK\cap H||K|.$$ It follows that $$|A|=|AK\cap H|\cdot|A\cap K|.$$ Since $K\le Z(G)$, $$C_{G}(A)=C_G(AK)=C_{G}((AK\cap H)\times K)=C_{H}(AK\cap H)\times K.$$ Hence
$$m_{G}(A)=|A|\cdot|C_{G}(A)|=|AK\cap H|\cdot|A\cap K|\cdot |C_{H}(AK\cap H)\times K|$$
$$=|AK\cap H|\cdot|C_{H}(AK\cap H)|\cdot|A\cap K|\cdot|K|=m_{H}(AK\cap H)\cdot m_{K}(A\cap K).$$
It follows that $|\mathrm{Im}(m_{G})|\leqslant|\mathrm{Im}(m_{H})|\cdot|\mathrm{Im}(m_{K})|$.\qed

\begin{lem}\label{rt}
Let $G=\langle a\rangle\rtimes \langle b\rangle$, where $o(a)=m,~o(b)=n$ and $(m,n)=1$. Then $|\mathrm{Im}(m_{G})|\leqslant\tau(m)\tau(n)-1$.
\end{lem}
\demo
Let $H\leq G$ and $|H|=kl$, where $k|m,~l|n$. By Lemma \ref{s1}, there exist $K\le H$ and $L\le H$ such that $|K|=k$ and $|L|=l$.
 Then $K\le \langle a\rangle$ by $\langle a\rangle \unlhd G$. Thus $K=\langle a^{\frac{m}{k}}\rangle$ and  $L=\langle (b^{g})^{\frac{n}{l}}\rangle$ for some $g\in G$ by Lemma \ref{s1}.
It follows that $$H=KL=\langle a^{\frac{m}{k}}, (b^{g})^{\frac{n}{l}}\rangle=\langle(a^{g})^{\frac{m}{k}}, (b^{g})^{\frac{n}{l}}\rangle=\langle a^{\frac{m}{k}},b^{\frac{n}{l}}\rangle^{g}.$$
By Lemma \ref{ge},
$$m_{G}(H)=m_{G}(\langle a^{\frac{m}{k}},b^{\frac{n}{l}}\rangle).$$ That is, the Chermak-Delgado measures of the same order subgroups of $G$ are equal. Since $m_{G}(G)=m_{G}(Z(G))$, $|\mathrm{Im}(m_{G})|\leqslant\tau(|G|)-1=\tau(m)\tau(n)-1$.\qed

\begin{lem}\label{e1}
Let $G\cong \mathrm{D}_{2n}=\langle a,b\mid a^{n}=1,b^2=1,[a,b]=a^{-2}\rangle$, where $n\geqslant 3$. Let $n=2^{l}\cdot k$, where $l\geqslant0$ and $2\nmid k$. Then

$(1)$ If $l=0$, then $|\mathrm{Im}(m_{G})|=2\tau(n)-1.$

$(2)$ If $l=1$, then $|\mathrm{Im}(m_{G})|=2\tau(n)-2.$

$(3)$ If $l\geqslant2$, then $|\mathrm{Im}(m_{G})|=2\tau(n)-4.$

\end{lem}
\demo
$(1)$ If $l=0$, then $|\mathrm{Im}(m_{G})|\leqslant 2\tau(n)-1$ by Lemma \ref{rt}. By calculating, $$m_{G}(1)=2n,~m_{G}(\langle a^{\frac{n}{s}}\rangle)=s\cdot |\langle a\rangle|=sn,$$ $$m_{G}(\langle b\rangle)=2\cdot|\langle b\rangle|=4,~m_{G}(\langle a^{\frac{n}{s}},b\rangle)=2s\cdot|\langle 1\rangle|=2s,$$ where $s\mid n$ and $s>1$. It is easy to see that $|\mathrm{Im}(m_{G})|=2\tau(n)-1.$

$(2)$ If $l=1$, then $Z(G)=\langle a^{\frac{n}{2}}\rangle$. Let $H\le G$ and $|H|=st$, where $s\mid 2^{2},~t\mid k$. If $s=1$, then $H=\langle a^{\frac{n}{t}}\rangle$. If $s=2$, then $H=\langle a^{\frac{n}{2t}}\rangle$ or $H=\langle a^{\frac{n}{t}}, a^{i}b\rangle$, where $i=0,1,\ldots,n-1$. If $s=4$, then $H=\langle a^{\frac{n}{2t}}, a^{i}b\rangle$. By calculating, we have
 \begin{center}
\begin{tabular}{|m{2.75cm}|m{1.95cm}|m{1.95cm}|m{1.95cm}|}
\hline
  $|H|$ & $H$ & $C_{G}(H)$ & $m_{G}(H)$\\\hline

  $1$ & $1$ & $\langle a,b\rangle$ & $2n$\\\hline

  $t(t>1)$ & $\langle a^{\frac{n}{t}}\rangle$ & $\langle a\rangle$ & $tn$\\\hline

  $2$ &  $\langle a^{i}b\rangle$ & $\langle a^{\frac{n}{2}}, a^{i}b\rangle$ & $8$\\\hline

  $2$ & $\langle a^{\frac{n}{2}}\rangle$ & $\langle a, b\rangle$ & $4n$\\\hline

  $2t(2<2t\leqslant 2k)$  & $\langle a^{\frac{n}{t}}, a^{i}b\rangle$ & $\langle a^{\frac{n}{2}}\rangle$ & $4t$\\\hline

  $2t(2<2t\leqslant 2k)$ & $\langle a^{\frac{n}{2t}}\rangle$ & $\langle a\rangle$ & $2tn$\\\hline

  $4$ & $\langle a^{\frac{n}{2}}, a^{i}b\rangle$ & $\langle a^{\frac{n}{2}}, a^{i}b\rangle$ & $16$ \\\hline

  $4t(t>1)$ & $\langle a^{\frac{n}{2t}}, a^{i}b\rangle$ & $\langle a^{\frac{n}{2}}\rangle$ & $8t$\\\hline

\end{tabular}.
\end{center}
Notice that, $$m_{G}(\langle a,b\rangle)=4n=m_{G}(\langle a^{\frac{n}{2}}\rangle)~\mbox{and}~m_{G}(\langle a^2, a^{i}b\rangle)=2n=m_{G}(1).$$
By calculating the number of possible values for $t$, and the number of distinct measures, we have $|\mathrm{Im}(m_{G})|=\tau(k)+2\tau(k)+\tau(k)-2=4\tau(k)-2=2\tau(n)-2$.

$(3)$ If $l\geqslant2$, then $Z(G)=\langle a^{\frac{n}{2}}\rangle$. Let $H\le G$ and $|H|=st$, where $s\mid 2^{l+1},~t\mid k$. If $s=1$, then $H=\langle a^{\frac{n}{t}}\rangle$. If $2<s\leqslant 2^{l}$, then $H=\langle a^{\frac{n}{st}}\rangle$ or $H=\langle a^{\frac{2n}{st}}, a^{i}b\rangle$, where $i=0,1,\ldots,n-1$. If $s=2^{l+1}$, then $H=\langle a^{\frac{n}{2^{l}t}}, a^{i}b\rangle$. By calculating, we have
 \begin{center}
\begin{tabular}{|m{2.75cm}|m{1.95cm}|m{1.95cm}|m{1.95cm}|}
\hline
  $|H|$ & $H$ & $C_{G}(H)$ & $m_{G}(H)$\\\hline

  $1$ & $1$ & $\langle a,b\rangle$ & $2n$\\\hline

  $t(t>1)$ & $\langle a^{\frac{n}{t}}\rangle$ & $\langle a\rangle$ & $tn$\\\hline

  $2$ &  $\langle a^{i}b\rangle$ & $\langle a^{\frac{n}{2}}, a^{i}b\rangle$ & $8$\\\hline

  $2$ & $\langle a^{\frac{n}{2}}\rangle$ & $\langle a, b\rangle$ & $4n$\\\hline

  $4$ &  $\langle a^{\frac{n}{2}}, a^{i}b\rangle$ & $\langle a^{\frac{n}{2}}, a^{i}b\rangle$ & $16$\\\hline

  $4$ & $\langle a^{\frac{n}{4}}\rangle$ & $\langle a\rangle$ & $4n$\\\hline

  $st(4<st\leqslant2^{l}k)$  & $\langle a^{\frac{2n}{st}}, a^{i}b\rangle$ & $\langle a^{\frac{n}{2}}\rangle$ & $2st$\\\hline

  $st(4<st\leqslant2^{l}k)$ & $\langle a^{\frac{n}{st}}\rangle$ & $\langle a\rangle$ & $stn$\\\hline

  $2^{l+1}t$ & $\langle a^{\frac{n}{2^{l}t}}, a^{i}b\rangle$ & $\langle a^{\frac{n}{2}}\rangle$ & $2^{l+2}t$\\\hline

\end{tabular}.
\end{center}
Notice that, $$m_{G}(\langle a^{\frac{n}{4}},a^{i}b\rangle)=16=m_{G}(\langle a^{\frac{n}{2}},a^{i}b\rangle),~m_{G}(\langle a,b\rangle)=4n=m_{G}(\langle a^{\frac{n}{2}}\rangle),$$
$$m_{G}(\langle a^2,a^{i}b\rangle)=2n=m_{G}(1),~m_{G}(\langle a^{\frac{n}{2}}\rangle)=4n=m_{G}(\langle a^{\frac{n}{4}}\rangle).$$
By calculating the number of possible values for $s$, $t$, and the number of distinct measures, we have $|\mathrm{Im}(m_{G})|=\tau(k)+2(\tau(2^{l})-1)\tau(k)+\tau(k)-4=2\tau(n)-4$.\qed

\begin{exam}
Let $n=60$ and $G$ be~a group~of~order $n$.

$(i)$ $\mathrm{Im}_{\rm max}(n,\mathcal{N})=12,~\mathrm{Im}_{\rm max}(n)=14;$

$(ii)$ If $G$ is non-solvable, then $G\cong A_{5}$ and $|\mathrm{Im}(m_{G})|=8;$

$(iii)$ $|\mathrm{Im}(m_{G})|=\mathrm{Im}_{\rm max}(n)$ if and only if $G\cong \mathrm{D}_{60}.$
\end{exam}

\demo
$(i)$ By Theorem B(1), we know that $\mathrm{Im}_{\rm max}(60,\mathcal{N})=12$.
A finite group $G$ is a non-nilpotent group of order $60$ if and only if $G$ is one of the
following pairwise non-isomorphic groups:

$(1)$ $A_{5}$;(non-solvable)

$(2)$ $\mathrm{C}_{15}\rtimes \mathrm{C}_{4}\cong\langle a,b\mid a^{15}=1,b^4=1,[a,b]=a\rangle;$

$(3)$ $\mathrm{C}_{15}\rtimes \mathrm{C}_{4}\cong\langle a,b\mid a^{15}=1,b^4=1,[a,b]=a^{-2}\rangle;$

$(4)$ $\mathrm{C}_{5}\times(\mathrm{C}_{3}\rtimes \mathrm{C}_{4})\cong\langle a,b,c\mid a^4=1,b^3=1,c^5=1,[b,a]=b,[a,c]=[c,b]=1\rangle;$

$(5)$ $\mathrm{C}_{3}\times(\mathrm{C}_{5}\rtimes \mathrm{C}_{4})\cong\langle a,b,c\mid a^4=1,b^3=1,c^5=1,[c,a]=c,[c,b]=[a,b]=1\rangle;$

$(6)$ $\mathrm{C}_{3}\times(\mathrm{C}_{5}\rtimes \mathrm{C}_{4})\cong\langle a,b,c\mid a^4=1,b^3=1,c^5=1,[c,a]=c^{-2},[c,b]=[a,b]=1\rangle;$

$(7)$ $\mathrm{C}_{5}\times A_{4}.$
$(8)$ $\mathrm{C}_{6}\times \mathrm{D}_{10};$
$(9)$ $\mathrm{C}_{10}\times S_{3};$
$(10)$ $\mathrm{D}_{60};$
$(11)$ $S_{3}\times \mathrm{D}_{10}.$

If $G\cong A_{5}$, then it is easy to verify that $A_{5}$ contains subgroups of orders $1$, $2$, $3$, $4$, $5$, $6$, $10$, $12$ and $60$, and that subgroups of the same order are conjugate. By Lemma \ref{ge}, the Chermak-Delgado measures of the same order subgroups of $G$ are equal. By calculating the Chermak-Delgado measures of subgroups of each order, we have $|\mathrm{Im}(m_{G})|=8$. Similarly, $|\mathrm{Im}(m_{A_{4}})|=4$.

If $G$ is one of the groups $(2)$--$(3)$, then $|\mathrm{Im}(m_{G})|\leqslant 11$ by Lemma \ref{rt}.

If $G$ is the group $(4)$, then  $|\mathrm{Im}(m_{G})|=|\mathrm{Im}(m_{\mathrm{C}_{5}})|\cdot|\mathrm{Im}(m_{\mathrm{C}_{3}\rtimes \mathrm{C}_{4}})|$ by Lemma \ref{lem1}(2). By Lemma \ref{a1}, $|\mathrm{Im}(m_{\mathrm{C}_{5}})|=2$. By Lemma \ref{rt}, $|\mathrm{Im}(m_{\mathrm{C}_{3}\rtimes \mathrm{C}_{4}})|\leqslant 5$. Then $|\mathrm{Im}(m_{G})|\leqslant 10$. Similarly, if $G$ is one of the groups $(5)$--$(6)$, then $|\mathrm{Im}(m_{G})|\leqslant 10$. If $G$ is the group $(7)$, then $|\mathrm{Im}(m_{G})|=8$.

If $G$ is the group $(8)$, then $|\mathrm{Im}(m_{G})|\leqslant|\mathrm{Im}(m_{\mathrm{C}_{6}})|\cdot|\mathrm{Im}(m_{\mathrm{D}_{10}})|$
by Lemma \ref{t2}. By Lemma \ref{a1} and Lemma \ref{lem1}(2), $|\mathrm{Im}(m_{\mathrm{C}_{6}})|=4$. Since $\mathrm{D}_{10}\cong \mathrm{C}_{5}\rtimes \mathrm{C}_{2}$, $|\mathrm{Im}(m_{\mathrm{D}_{10}})|\leqslant 3$ by Lemma \ref{rt}. Thus $|\mathrm{Im}(m_{G})|\leqslant 12$ . Similarly, if $G$ is the group $(9)$, then $|\mathrm{Im}(m_{G})|\leqslant 12$.

If $G$ is the group $(10)$, then $|\mathrm{Im}(m_{G})|=14$ by Lemma \ref{e1}(2).

If $G$ is the group $(11)$, then let $S_{3}\times \mathrm{D}_{10}=\langle a,b,c,d\mid a^{3}=1,b^2=1,c^5=1,d^2=1,[a,b]=a^{-2},[c,d]=c^{-2},[a,c]=[b,c]=[a,d]=[b,d]=1\rangle$. By calculating, we have
 \begin{center}
\begin{tabular}{|m{0.75cm}|m{1.5cm}|m{1.5cm}|m{1.25cm}|m{0.75cm}|m{1.5cm}|m{1.5cm}|m{1.25cm}|}
\hline
  $|H|$ & $H$ & $C_{G}(H)$ & $m_{G}(H)$ & $|H|$ & $H$ & $C_{G}(H)$ & $m_{G}(H)$\\\hline
  $1$ & $1$ & $G$& 60 & $10$ & $\langle c,d\rangle$ & $\langle a,b\rangle$ & $60$\\\hline
  $2$ & $\langle a^{i}b\rangle$ & $\langle a^{i}b,c,d\rangle$ & $40$ & $10$ & $\langle c,da^{i}b\rangle$ & $\langle a^{i}b\rangle$ & $20$\\\hline
  $2$ & $\langle c^{j}d\rangle$ & $\langle a,b,c^{j}d\rangle$ & $24$ & $10$ & $\langle c,a^{i}b\rangle$ & $\langle c,a^{i}b\rangle$ & $100$\\\hline

  $2$ & $\langle a^{i}bc^{j}d\rangle$ & $\langle a^{i}b,c^{j}d\rangle$ & $8$ & $15$ & $\langle a,c\rangle$ & $\langle a,c\rangle$ & $225$\\\hline

  $3$ & $\langle a\rangle$ & $\langle a,c,d\rangle$ & $90$ & $12$ & $\langle a,b,c^{j}d\rangle$ & $\langle c^{j}d\rangle$ & $24$\\\hline

  $5$ & $\langle c\rangle$ &  $\langle a,b,c\rangle$  & $150$ & $20$ & $\langle a^{i}b,c,d\rangle$ & $\langle a^{i}b\rangle$ & $40$\\\hline

  $4$ & $\langle a^{i}b,c^{j}d\rangle$ &  $\langle a^{i}b,c^{j}d\rangle$ & $16$ & $30$ & $\langle a,c,d\rangle$ & $\langle a\rangle$ & $90$\\\hline

  $6$ & $\langle a,b\rangle$ &  $\langle c,d\rangle$ & $60$ &  $30$ & $\langle a,b,c\rangle$ & $\langle c\rangle$ & $150$\\\hline

  $6$ & $\langle a,c^{j}d\rangle$ &  $\langle a,c^{j}d\rangle$ & $36$ &  $30$ & $\langle a,c,bd\rangle$ & $1$ & $30$ \\\hline

  $6$ & $\langle a,bc^{j}d\rangle$ &  $\langle c^{j}d\rangle$ & $12$ & $60$ & $\langle a,b,c,d\rangle$ & $1$ & $60$\\\hline
\end{tabular}
\end{center}
where $i=0,1,2$ and $j=0,1,2,3,4$. Hence $|\mathrm{Im}(m_{G})|=13$.

By comparison, we can see that $\mathrm{Im}_{\rm max}(60)=14$.

$(ii)-(iii)$ The deduction can be directly obtained from the proof process of $(i)$.\qed

\bigskip

According to the above example, we indicate some natural open problems concerning the above study.

\begin{prob}
Let $n$ be a positive integer.

$(1)$ What is the value of $\mathrm{Im}_{\rm max}(n)?$

$(2)$ Determine the finite groups $G$ with $|\mathrm{Im}(m_{G})|=\mathrm{Im}_{\rm max}(n).$

$(3)$ Are the groups in $(2)$ are solvable$?$
\end{prob}

\end{document}